\documentclass[conference]{IEEEtran}
\usepackage{graphics} 
\usepackage{epsfig} 
\usepackage{amsmath} 
\usepackage{amssymb}  
\usepackage[caption=false,font=footnotesize]{subfig}
\usepackage[hidelinks]{hyperref} 

\ifCLASSINFOpdf
\else
\fi
\begin{document}
%
\title{Attack-Aware Multi-Sensor Integration Algorithm\\for Autonomous Vehicle Navigation Systems}


\author{
\IEEEauthorblockN{Sangjun Lee}
\IEEEauthorblockA{Department of Computer and \\Information Technology\\
Purdue University\\
West Lafayette, IN 47907\\
Email: lee1424@purdue.edu}
\and
\IEEEauthorblockN{Yongbum Cho}
\IEEEauthorblockA{School of Mechanical Engineering\\
Purdue University\\
West Lafayette, IN 47907\\
Email: cho148@purdue.edu}
\and
\IEEEauthorblockN{Byung-Cheol Min}
\IEEEauthorblockA{Department of Computer and \\Information Technology\\
Purdue University\\
West Lafayette, IN 47907\\
Email: minb@purdue.edu}
}


%


\maketitle

\begin{abstract}
In this paper, we propose a fault detection and isolation based attack-aware multi-sensor integration algorithm for the detection of cyberattacks in autonomous vehicle navigation systems. The proposed algorithm uses an extended Kalman filter to construct robust residuals in the presence of noise, and then uses a parametric statistical tool to identify cyberattacks. The parametric statistical tool is based on the residuals constructed by the measurement history rather than one measurement at a time in the properties of discrete-time signals and dynamic systems. This approach allows the proposed multi-sensor integration algorithm to provide quick detection and low false alarm rates for applications in dynamic systems. An example of INS/GNSS integration of autonomous navigation systems is presented to validate the proposed algorithm by using a software-in-the-loop simulation.
\end{abstract}


%
\IEEEpeerreviewmaketitle

\section{Introduction}
Security of Cyber-Physical Systems (CPS) has garnered significant attention as a major issue with regard to autonomous vehicles. Today's autonomous vehicles enable the deployment of safety technologies, such as automatic emergency braking, collision warning, and Vehicle-to-Everything technologies. In the near future, these systems will be available in all vehicles to help achieve zero fatalities, zero injuries, and zero accidents. However, behind the great potential of these innovations, a new challenge of ensuring security from cyberattacks needs to be addressed.

A typical autonomous vehicle receives and transmits a great deal of information between sensors, actuators, and the electronic control units, all providing access for attackers \cite{petit2015potential}. From this point of view, cybersecurity is imperative. Units that govern safety should be protected from malicious attacks, unauthorized access, or dubious activities, all of which could cause harmful outcomes. For example, an autonomous vehicle's navigation system must be secured because it controls real-time position data directly linked to the physical behavior of the vehicle. We have a real-world example \cite{miller2015remote} in which a hack was able to remotely hijack a car, and other examples \cite{kerns2014unmanned,shepard2012evaluation} in which unmanned aerial vehicles were captured and controlled via Global Positioning System (GPS) signal spoofing. Practical studies on the analysis of security vulnerabilities of autonomous vehicles have been discussed in \cite{miller2014survey,amoozadeh2015security}. Similarly, an extensive study of potential cybersecurity threats to autonomous vehicles was published in the open literature \cite{petit2015potential}. This study presented many possible attack methods and identified that sensor spoofing and false data injection could result in the worst safety related issue.

Securing autonomous vehicles' safety is challenging because it requires the full knowledge of applications that consist of numerous hardware and multi-layered architectures \cite{urbina2016survey}. For instance, an autonomous vehicle navigation system is generally comprised of multiple sensors such as Inertial Navigation System (INS) and Global Navigation Satellite System (GNSS). These two different types of sensors have inherent limitations so that integration methodologies for such systems have been widely introduced to combine the advantages of both technologies \cite{rezaei2007kalman}. However, an integrated system does not have any safety functions against cyberattacks, leaving it highly vulnerable. Additionally, the lack of knowledge of multi-sensor integration makes autonomous vehicles more exposed to cyberattacks. A fault tolerant multi-sensor perception system was presented to provide fault-free inputs for critical functions of mobile robots \cite{bader2017fault}. All of the previously mentioned studies suggest that there are rapidly growing needs for ensuring cybersecurity in autonomous vehicles.

One of the common approaches for achieving security guarantee is the Fault Detection and Isolation (FDI) method. This approach has been widely studied in various applications such as spacecraft \cite{pirmoradi2007efficient}, aircraft \cite{abbaspour2016detection}, power system \cite{mohanty2008cumulative}, and automobile \cite{huang2015design}. In general, a fault detection algorithm generates a residual and compares it with a predefined threshold. If the residual exceeds the threshold, the algorithm reveals a fault and an alarm is triggered. In this manner, abnormal dynamic behavior and abrupt system changes caused by cyberattacks can be detected. The authors in \cite{hwang2010survey,patcha2007overview} have presented a remarkable comparison of existing residual generation algorithms and threshold determination techniques.

The primary focus of attack detection for dynamic systems is to generate residuals and design decision rules based upon these residuals. Ideal residuals would be zero under normal operation when there is no attack. However, residuals are subject to the presence of noise and unknown errors in real-world applications \cite{gustafsson2000adaptive}. For this reason, it is challenging to generate robust residuals that are insensitive to noise and uncertainties yet sensitive to attacks in order to provoke a quick alarm \cite{basseville1993detection}. Optimal filters and state observers have been proposed to generate a sequence of residuals that resemble white noise in normal operation \cite{oonk2014extended,marino2015discrete}. After residual generation, an attack alarm will be triggered at the moment residuals exceed the threshold. Another challenge here is to determine the threshold limit. This is a fundamental limitation of attack detection because determining thresholds is a compromise between detecting true attacks and avoiding false alarms. Some studies have proposed statistical approaches to generate an adaptive threshold in order to avoid false alarms \cite{fillatre2014statistical,pradhan2006moving}. Others have used a hypothesis test with Boolean questions to determine system attacks \cite{murguia2016characterization}.

Although the aforementioned studies have presented various strategies and solutions for attack detection, there are still questions to address. The lack of knowledge of interaction among sensors, actuators, and electronic control units increases the possibilities of being compromised by unidentified source. Therefore, the following research questions can be raised:
\begin{itemize}
\item {How will the driver know when he or she has to take back control from full self-driving mode due to security breach?}
\item {How will the system identify possible attacks against multi-sensors that are tightly coupled instead of a single sensor?}
\item {How will the system present state estimates as close to the true value as possible in the presence of noise without compromising response time or sensitivity?}
\end{itemize}
To provide answers to the questions, this paper focuses on possible attacks on the autonomous vehicle navigation systems. It is a highly vulnerable system because it handles signals from external sources. Thus, this study determines that a vehicle's navigation system is being attacked if any abrupt change or unexpected dynamic behavior has been identified by a proposed algorithm. We assume that system alterations are caused by false data injection attacks, corrupted signal reading, sensor failure, or any combination of these.

To summarize, the main contributions of this work are as follows:
	\begin{enumerate}
	\item Development of an attack-aware multi-sensor integration algorithm for the autonomous vehicle navigation system;
	\item Generation of robust residuals in the presence of uncertainties;
	\item Design of a parametric statistical test that enables the proposed algorithm to generate a quick detection alarm and low false alarm rate;
	\item Application of the proposed algorithm to the detection of attack on INS/GNSS integration of autonomous vehicles;
	\item Verification of the application in a customized software-in-the-loop simulation.
	\end{enumerate}

The rest of this paper is organized as follows. In Section \ref{sec:met}, an attack-aware multi-sensor integration is developed with the strategies of residual generation and threshold determination. In Section \ref{sec:application}, the proposed attack detection algorithm with an application to the autonomous navigation system is introduced and a simulation is designed to validate it. Finally, conclusions and future works are discussed in Section \ref{sec:con}.

\section{Problem Formulation}\label{sec:met}
\begin{figure}[t]
\centering
\includegraphics[width=3.5 in]{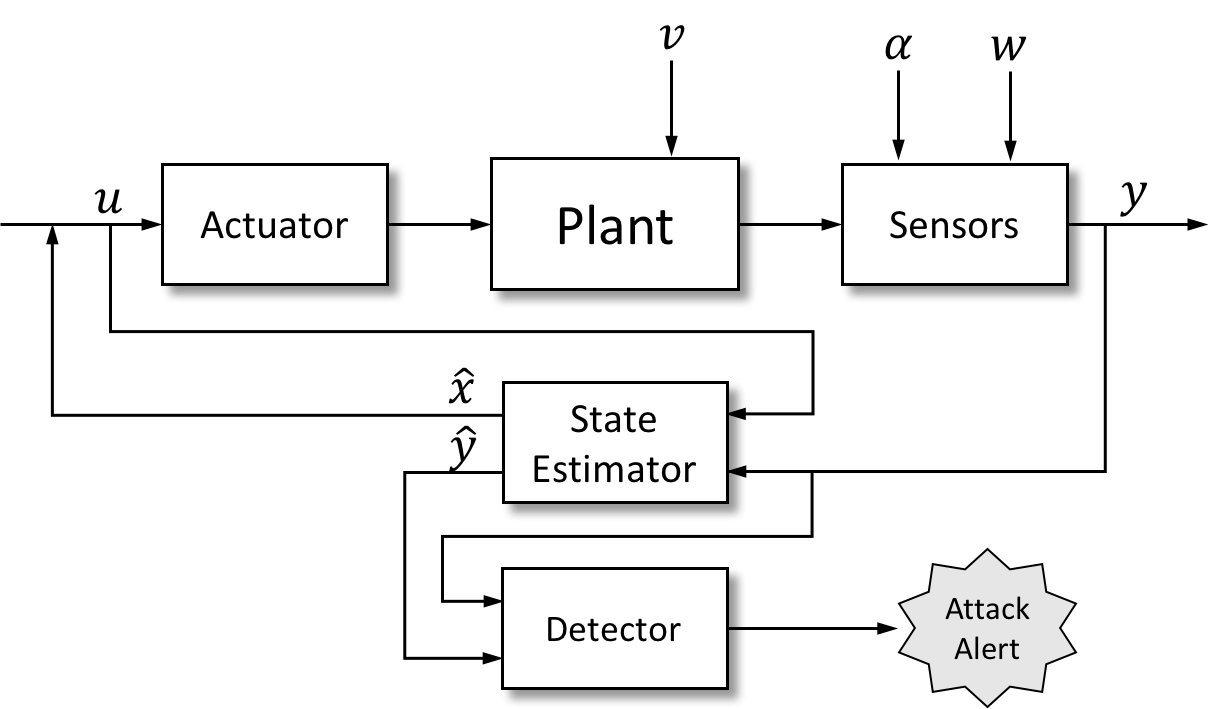}
\caption{An overview of the proposed attack-aware multi-sensor integration system. An attack is introduced to the sensor.} 
\label{fig:blocks}
\end{figure}

This section provides a Kalman filter-based estimation for a multi-sensor integration and detection algorithm. The system model that we consider is illustrated in Fig. \ref{fig:blocks}. 
The actuator sends a command to the plant in accordance with the control input and then the sensors measure some of the states. These states are fed into the state estimator to predict the states. Lastly, the detector determines if there is an attack on the sensor through comparison between state estimations and sensor measurements.

\subsection{Attack Model} 
We investigate attacks in the state or measurement equation of a discrete liner time-invariant (LTI) system represented by a state-space model. The state-space model with given matrices $A, B$, and $C$ is given as
	\begin{align}
	x(k+1) & = Ax(k) + Bu(k) + \nu(k) \label{eq:ss_normal_sys} \\
	y(k) &= Cx(k) + \omega(k), \label{eq:ss_normal_mea}
	\end{align}
where $x \in \mathcal{R}^{n}$, $y \in \mathcal{R}^{m}$, and $u \in \mathcal{R}^{r}$ represent state vector, output vector, and control input vector, respectively, and where $\nu$ and $\omega$ are process and measurement noise that are represented by two independent white noise sequences with covariance matrices $Q$ and $R$, respectively.
If a sensor is being compromised that means unknown signals have been injected, added, or modified to the sensor, the LTI system (\ref{eq:ss_normal_sys}) and (\ref{eq:ss_normal_mea}) can be written as follows:
	\begin{align}\label{eq:ss_attack}
    \begin{split}
	x(k+1) & = Ax(k) + Bu(k) + \nu(k) \\
	y_{\alpha}(k) &= Cx(k) + \alpha(k) + \omega(k),
    \end{split}
    \end{align}
where $\alpha \in \mathcal{R}^{m}$ denotes additive attacks on a sensor and the state with the subscript $\alpha$ represents the system after an attack occurs. The key idea behind this is that the difference induced by attacks would be observable from the detection algorithm in the presence of uncertainties.

\subsection{Multi-sensor Integration}\label{sec:met:integration}
\begin{figure}[t]
\centering
\includegraphics[width=3.4 in]{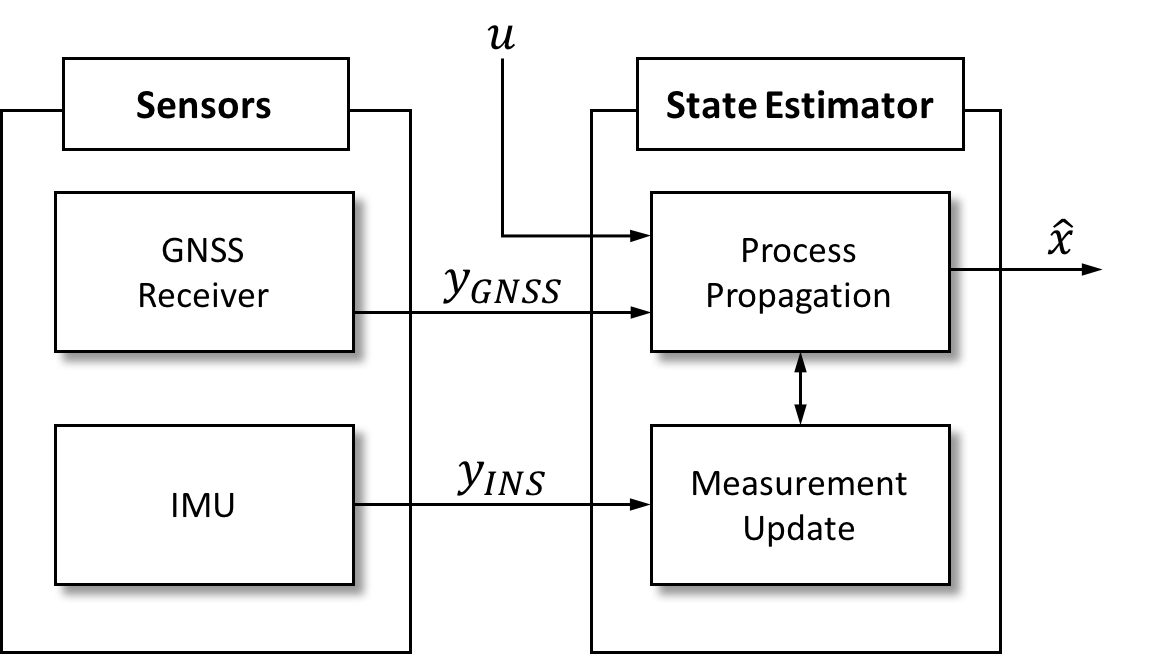}
\caption{Subsystems of the sensor and the state estimator. These subsystem are used in the Kalman filter-based multi-sensor integration.}
\label{fig:integration}
\end{figure}

A state estimator is designed to predict states from available measurements since not all the states of a system are observable in real-world applications. Two typical navigation solutions of autonomous vehicles, INS and GNSS measurements, are considered as shown in Fig. \ref{fig:integration}. An INS uses an Inertial Measurement Unit (IMU) to track the position, velocity, and orientation of a vehicle relative to an initial point, orientation, and velocity. A GNSS provides satellite signals that can be processed in a GNSS receiver, allowing the receiver to estimate its current position and velocity. The advantages of both technologies can be combined by fusing these navigation solutions. There are no states directly affected by the INS measurements or the GNSS measurements in the system model (\ref{eq:ss_normal_sys}), but they interact through the output vector (\ref{eq:ss_normal_mea}) determined by the measurement models:
	\begin{equation}\label{eq:output}
	y =
	\begin{bmatrix}
	y_{\text{GNSS}} \\ y_{\text{INS}}
	\end{bmatrix}.
	\end{equation}

Under the assumption that the system will stay in the steady-state until any attacks happen, it enables the system to identify any abrupt changes on sensor measurements. An estimator dynamics given by the following steady-state Kalman filter is considered:
	\begin{equation}
    \hat{x}(k+1) = A\hat{x}(k) + Bu(k) + K[y(k) - \hat{y}(k)],
    \end{equation}
where Kalman gain is $K = PC^{T}(CPC^{T} + R)^{-1}$ with the covariance matrix given by $P = A[P-PC^{T}(CPC^{T}+R)^{-1}CP]A^{T} + Q$. Note that the detectability of $(A,C)$ ensures the existence of such estimator.
This multi-sensor integration gives a continuous position estimation and achieves precise vehicle control.
	
\subsection{Detection Algorithm}\label{sec:met:detection}
\begin{figure}[t]
\centering
\includegraphics[width=3.4 in]{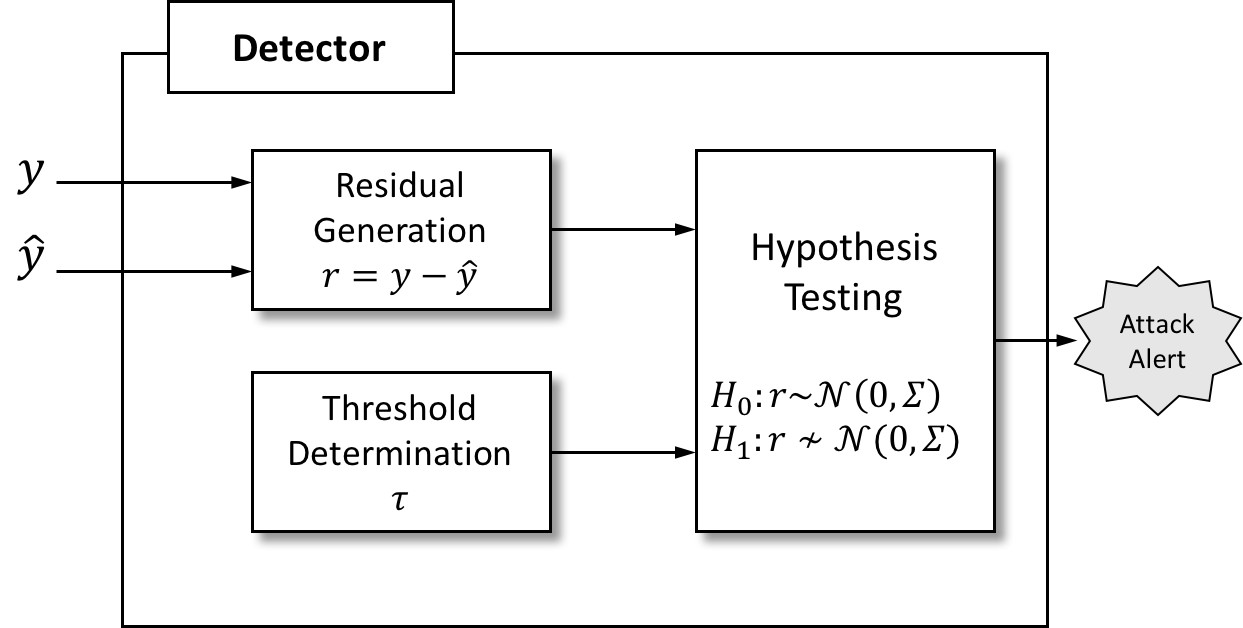}
\caption{A subsystem of the detector. A hypothesis testing determines the system functionality.}
\label{fig:detection}
\end{figure}
The main idea of the detection capability is to generate robust residuals to uncertainties and determine sensitive thresholds to false alarm. As shown in Fig. \ref{fig:detection}, the detector determines the system condition at each time step through statistical hypothesis testing that compares the residual and threshold generated. The residual is the difference between the actual measurements and the estimates. A sequence of the residuals is defined as
	\begin{equation}\label{eq:resi}
    r(k) = y_{\alpha}(k) - \hat{y}(k).
    \end{equation}

The residuals evolve with the output estimate given by $\hat{y}(k) = C \hat{x}$ and the estimation error defined as $e(k) = x- \hat{x}$. The residual dynamics is written as
	\begin{equation}
    r(k+1) =Ce(k+1) + \alpha(k+1),
    \end{equation}
where the estimation error dynamics given by $e(k+1) = (A-KC) e(k)$. Regardless of the availability of prior information, the residual is ideally zero before the attack and nonzero after the attack. Thus, if the system is under normal operation, the mean of the residuals will be zero and the covariance will have a value:
	\begin{align}
    E[r(k+1)]&=0 \\
    \Sigma[r(k+1)]&=CPC^{T} + R,
    \end{align}
where $E[\cdot]$ denotes the expected value and $\Sigma[\cdot]$ denotes the covariance matrix. The system is able to construct a two-sided hypothesis testing to make a decision at each time step when given a set of samples. It determines the system's abnormal behavior with the null hypothesis of normal operation and the alternative hypothesis of abnormal operation as follows:
	\begin{align}\label{eq:htest}
    \begin{split}
    \mathcal{H}_{0}: r(k) & \sim \mathcal{N}(0,\Sigma) \\
    \mathcal{H}_{1}: r(k) & \nsim \mathcal{N}(0,\Sigma),
    \end{split}
    \end{align}
where $\mathcal{N}(\sigma,\Sigma)$ denotes the probability density function of the Gaussian random variable with mean $\sigma$ and covariance matrix $\Sigma$. The test will continue as long as the decision favors the hypothesis $\mathcal{H}_{0}$ while the test will be stopped and restarted if the decision favors $\mathcal{H}_{1}$. Decision rules for rejecting the null hypothesis are based on the Cumulative Summation (CUSUM) algorithm which was introduced by Page \cite{page1954continuous}. In case of the system described in (\ref{eq:htest}), the two-sided CUSUM test is defined as
	\begin{align}\label{eq:cusum}
	S(k+1)=
	\begin{cases}
	\max{(0,S(k) + |r(k+1)|)} & \text{ if } S(k) \le \tau(k) \\
	0 \text{ and } k_{\alpha}=k & \text{ if } S(k) > \tau(k). 
	\end{cases}
	\end{align}
    
The null hypothesis is rejected if the test statistics $S$ is greater than the threshold $\tau$. In this case, the test provides an attack alarm time $k_{\alpha}$ and the test starts over. The null hypothesis is accepted if the test statistics $S$ is less than or equal to the threshold $\tau$. The test continues without stopping in this case. In practice, this test collects a number of samples and calculates their weighted sum to detect a significant change in the mean of samples.
Note that a selection of the sample size $N = 1,2, \cdots, k+1$ is to find a balance between response time and sensitivity while a selection of the threshold is to find a balance between sensitivity and a false alarm rate.

\section{Application to Navigation System of Autonomous Vehicles}\label{sec:application}
In this section, the proposed attack-aware integration algorithm is applied to a navigation system of an autonomous vehicle in the presence of uncertainties and unknown attacks on sensors. It is imperative that units such as the navigation system that govern safety are protected from malicious attacks, unauthorized access, or dubious activities. This is because a small change could result in significant changes in behavior. For the simulation studies, a vehicle model and sensor models are considered. An EKF is used for online estimation and multi-sensor integration as described in Section \ref{sec:met:integration}. According to the detection algorithm in Section \ref{sec:met:detection}, a significant change in the mean is detected and indicates an attack. A numerical simulation with a robotic simulator demonstrates the performance of the proposed algorithm. The following assumptions are considered through the simulation: 
no attack on multiple sensors at a time;
a random attack injection time;
an arbitrary magnitude of attack but greater than sensor biases.

\subsection{Design of Software-in-the-loop Simulation (SILS)}
\begin{figure}[t]
\centering
\subfloat[Feedback loop enclosing ROS environment variables. An attack-aware multi-sensor integration algorithm is built in the MATLAB node, and a robotic simulator runs on the Gazebo node. Each node is able to exchange data via topic subscription and publication with unique types of ROS messages.]{\includegraphics[width=3.2 in]{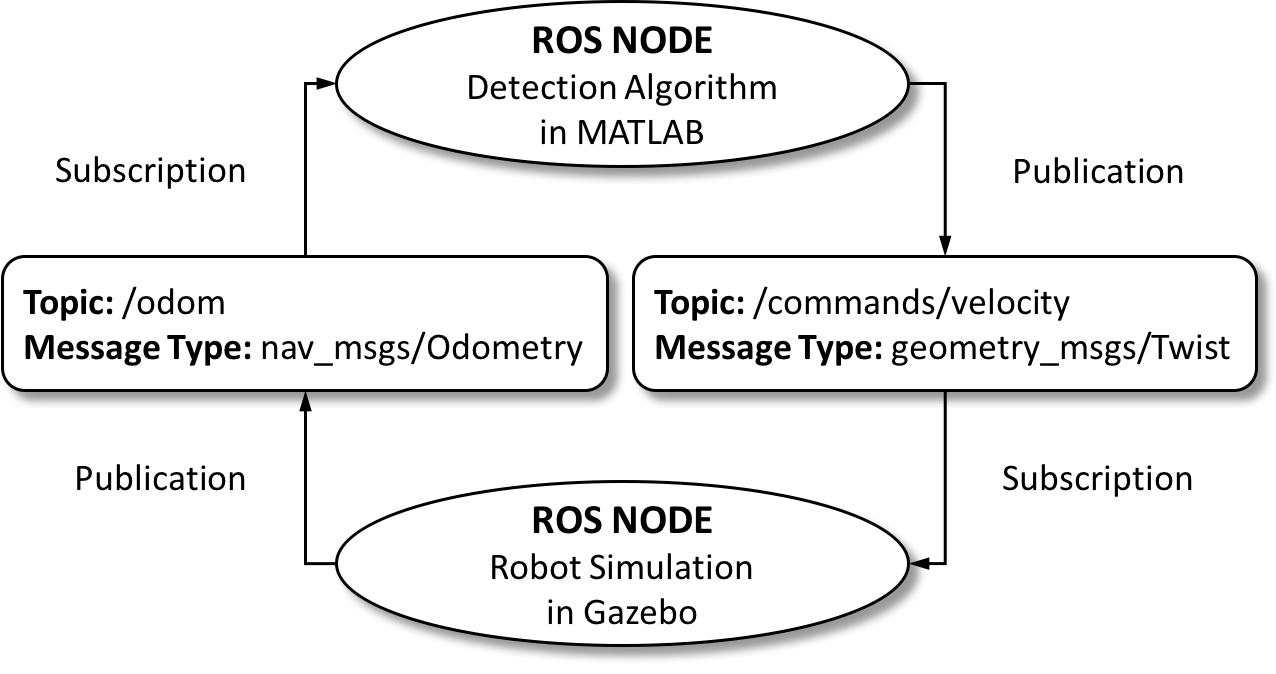}
\label{fig:rosnodes}}
\hfil
\subfloat[Gazebo simulation environment. A vehicle follows the desired path which is a straight line from the initial location at the bottom left to the home at the top right.]{\includegraphics[width=3.2 in]{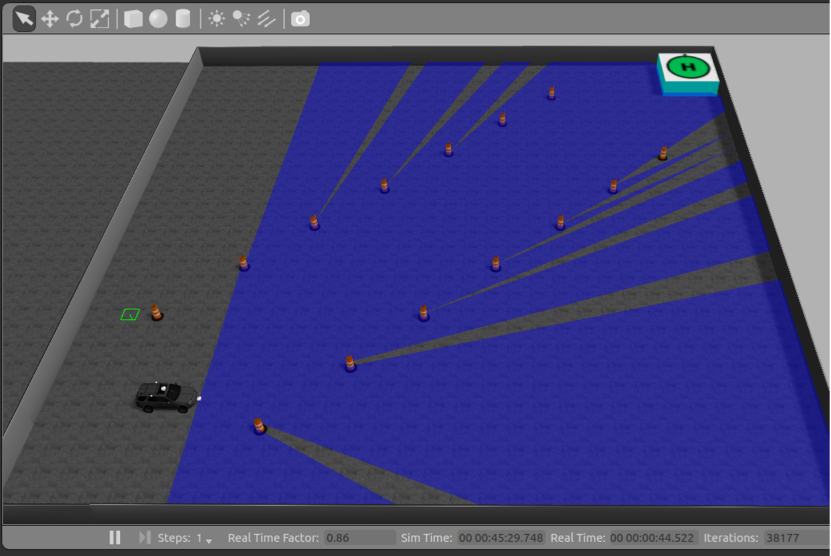}
\label{fig:gazebo}}
\caption{Software-in-the-loop simulation environment.}
\label{fig:sils}
\end{figure}

A software-in-the-loop simulation is designed to evaluate the proposed algorithm with an application of autonomous vehicles. The complete model of the simulation is illustrated in Fig. \ref{fig:sils}. The simulation runs on Robot Operating System (ROS), and it includes two ROS nodes as shown in Fig. \ref{fig:rosnodes}. One node is MATLAB that runs the multi-sensor integration and the detection algorithm, and another node is Gazebo that runs the robotic simulator in a customized world as shown in Fig. \ref{fig:gazebo}. Each node is able to create a unique topic in ROS message type. It enables each node to exchange data via topic subscription and publication without conflict.

For the model of an autonomous vehicle in the simulation, the CAT Vehicle, a full-sized model of Ford Escape developed by the Compositional Systems Laboratory at the University of Arizona \cite{catvehicle}, was used. It was actuated to be controllable through unique ROS topics. The simulation started with providing a set of desired waypoints to the mathematical model of the vehicle in MATLAB. The model then published the velocity commands subscribed by the robotic simulator in Gazebo. The CAT Vehicle in Gazebo followed the commands and published its local position data subscribed by the position controller in MATLAB to generate a new velocity command for the next time step. This feedback loop ran continuously and recursively until the vehicle reached the final destination regardless of attacks, and the sampling rate was 10 Hz.

\subsection{Implementation}
A loosely coupled INS/GNSS navigation model with a vehicle model is considered to represent an autonomous vehicle navigation system. Firstly, an EKF-based multi-sensor integration is developed for the residual generation. It is comprised of the state model and the measurement model.
Consider the equation of motion for the vehicle is governed by the following dynamics:  
	\begin{align}\label{eq:con}
    \begin{split}
	\dot{x} &= v_{x} \cos{\theta} - v_{y} \sin{\theta} \\
	\dot{y} &= v_{x} \sin{\theta} + v_{y} \cos{\theta},
    \end{split}
	\end{align}
where $x, y, v_{x}$, and $v_{y}$ represent the position along the eastern axis, the position along the northern axis, the velocity along the eastern axis, and the velocity along the northern axis, respectively.
The yaw angle is represented as $\theta$. The continuous time state equations can be discretized with the sampling time $T$ which gives the nonlinear discrete-time state model under normal operation as:

	\begin{align}
    \begin{split}
    x(k+1) &= x(k) + T v_{x}(k) \cos{\theta(k)} - T v_{y}(k) \sin{\theta(k)} \\
    y(k+1) &= y(k) + T v_{x}(k) \sin{\theta(k)} + T v_{y}(k) \cos{\theta(k)} \\
    \theta(k+1) &= \theta(k) + T \dot{\theta}(k) \\ 
    v_{x}(k+1) &= v_{x}(k) + T a_{x}(k) \\
    v_{y}(k+1) &= v_{y}(k) + T a_{y}(k) \\
    \dot{\theta}(k+1) &= \dot{\theta}(k) \\
    a_{x}(k+1) &= a_{x}(k) \\
    a_{y}(k+1) &= a_{y}(k) \\
    b_{\dot{\theta}}(k+1) &= b_{\dot{\theta}}(k) \\
    b_{a_{x}}(k+1) &= b_{a_{x}}(k) \\
    b_{a_{y}}(k+1) &= b_{a_{y}}(k),
    \end{split}
    \end{align}
and the linear measurement model under normal operation is given by
	\begin{align}\label{eq:mea_normal}
    \begin{split}
    y_{x}(k+1) &= x(k) \\
    y_{y}(k+1) &= y(k) \\
    y_{\theta}(k+1) &= \theta(k) \\
    y_{\dot{\theta}}(k+1) &= \dot{\theta}(k) + b_{\dot{\theta}}(k) \\
    y_{a_{x}}(k+1) &= a_{x}(k) + b_{a_{x}}(k) \\
    y_{a_{y}}(k+1) &= a_{y}(k) + b_{a_{y}}(k),
    \end{split}
    \end{align}
where $a$ and $b$ represent the acceleration and bias, respectively. Note that the process noise $\nu$ and measurement noise $\omega$ are additive to each equation. These models are linearized to correspond with the state-space model in  (\ref{eq:ss_normal_sys}) and  (\ref{eq:ss_normal_mea}) by using the state and measurement Jacobian matrices. In addition, initial states $x(0)$, state error covariance $P$, process noise covariance $Q$, and measurement noise covariance $R$ are carefully chosen according to hardware specifications.
The models in (\ref{eq:con})-(\ref{eq:mea_normal}) integrate multiple sensors to predict the vehicle states under normal operation. This integrated architecture ensures that a continuous navigation solution is always produced, regardless of the existence of attacks. Following the state estimation under normal condition, the system under attack  (\ref{eq:ss_attack}) is considered. These two different measurement models are used for the residual generation in (\ref{eq:resi}). The decision rules in (\ref{eq:cusum}) then determine if there is a significant change in the vehicle position at each time step. It is verified in the following section.

\subsection{Results}
\begin{figure}[t]
\centering
\includegraphics[width=3.5in, height=2.5in]{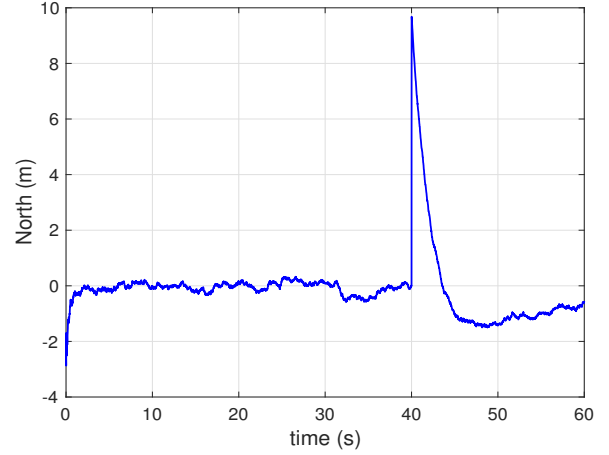}
\caption{North position estimation error corresponding to an attack in the vehicle navigation system. A peak is observed around 40 seconds but it does not indicate that the peak has been caused by the attack.}
\label{fig:EstErr}
\end{figure}
\begin{figure}[t]
\centering
\includegraphics[width=3.5in, height=2.5in]{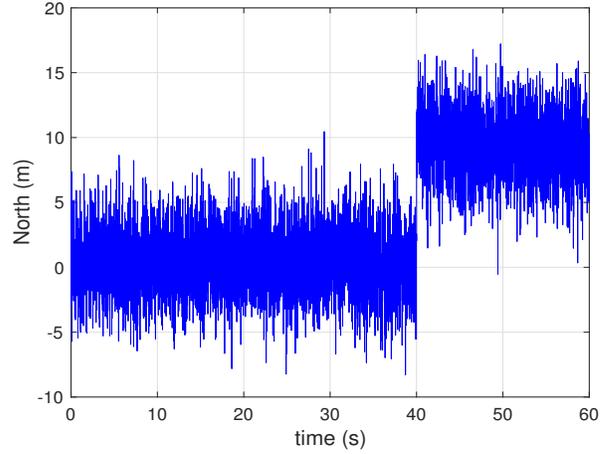}
\caption{North position measurement error corresponding to an attack in the vehicle navigation system. The measurement error jumped around the 40 second mark by approximately 10 meters but it does not guarantee that the shift occurred due to the attack.}
\label{fig:MeaErr}
\end{figure}
\begin{figure}[t]
\centering
\includegraphics[width=3.5in, height=2.5in]{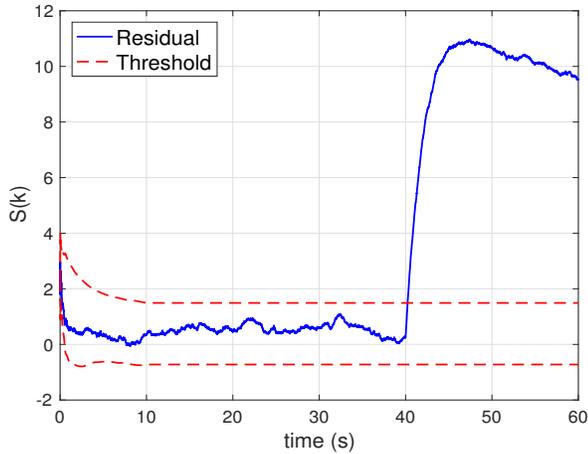}
\caption{Test statistics evolution corresponding to an attack in the vehicle navigation system. The proposed algorithm identified a significant change of the residuals that exceeds the upper limit of the threshold as soon as the attack was initiated at 40 seconds.}
\label{fig:ResEvol}
\end{figure}

During the simulation, an attack was introduced at the GNSS receiver at 40 seconds to test if the proposed detection algorithm can identify the attack. A separate function from the detection algorithm injected the attack into the receiver measurement if the simulation clock reached 40 seconds, and there was no data exchange with the detection algorithm. The magnitude of the attack was 10 meters, which is larger than the GNSS receiver bias.

The estimation error in Fig. \ref{fig:EstErr} shows the estimation performance of the multi-sensor integration. There are quite small errors, which means it provides a continuous and high-bandwidth navigation solution, until a peak around 40 seconds. The peak may imply that there was an attack around 40 seconds but it is insufficient evidence to determine that the peak was due to an attack. This is because an attack is not the only cause of a peak during state estimation. For example, it can be caused by signal attenuation, data loss, time delay, bursty packet dropping, etc. 
Similarly, the measurement error in Fig. \ref{fig:MeaErr} indicates that there was an abrupt shift around 40 seconds on the north sensor measurement. This is not sufficient to determine if an attack was introduced because it is unable to verify where the shift originates. Consequently, one can indicate a suspicious jump or shift from the multi-sensor integration but it is insufficient to determine that there is an attack on the vehicle. 
On the other hand, the evolution of the test statistics in Fig. \ref{fig:ResEvol} clearly shows that there was a significant change that caused the residual to jump the upper bound of the threshold around the 40 second mark. The test statistics were calculated by (\ref{eq:cusum}), and the upper and lower bounds of the threshold were generated by using the weighted sum of the first 10 samples. Based upon these parameters, the detector in the navigation system determined that there was an attack around 40 seconds when the residual went above the upper limit of the threshold, and the corresponding time was automatically generated. It was 40.2 seconds in this simulation, two time steps behind the attack (i.e. an attack was injected at $k=400$ but $k_{\alpha}=402$), a fairly quick detection because it was only two sampling steps behind the actual attack. In addition, there were a number of ups and downs prior to the attack but they stayed within the threshold boundary, allowing the detection algorithm to avoid a false alarm. Thus in this application, using the proposed attack-aware multi-sensor integration system provides a method to detect an attack as quickly as possible with no false alarm.

\section{Conclusion} \label{sec:con}
This research presented a statistical approach to the problem of attack detection on the multi-sensor integration of autonomous vehicle navigation systems. Starting with a state-space model of the system under attack, a parametric statistical tool with a multi-sensor integration strategy was developed to identify an attack. Finally, a simulation was designed to verify the proposed detection system and results were presented. A few limitations in this study remain: 1) the detection system was unable to identify an attack that was smaller than the sensor bias, but the vehicle was still under the control, and 2) the detection system was unable to detect an attack if any change occurred at the very beginning of samples. These remaining research questions will be addressed in the future.

\bibliographystyle{IEEEtran}
\bibliography{references}

\end{document}